\documentclass[10pt]{article}
\usepackage{amssymb}
\usepackage{amsmath,amsthm}
\usepackage{mathtools}
\usepackage{bm}
\usepackage{extarrows}
\usepackage{indentfirst}
\usepackage[english]{babel}
\usepackage{tikz}
\usepackage[numbers]{natbib}
\usepackage[utf8]{inputenc}
\makeatletter

\newcommand{\Rmnum}[1]{\expandafter\@slowromancap\romannumeral #1@}
\makeatother

\def\span{\text{span}}
\def\basis{\text{basis}}

\newtheorem{thm}{Theorem}[section]

\newtheorem{lem}[thm]{Lemma}
\newtheorem{prop}[thm]{Proposition}

\newtheorem{prob}[thm]{Problem}
\newtheorem{claim}{Claim}

\begin{document}
\title{Maximal fractional cross-intersecting families~\thanks{The work was supported by National Natural Science Foundation of China (No. 12071453) and the National Key R and D Program of China(2020YFA0713100).}}
\author{Hongkui Wang$^a$, \quad Xinmin Hou$^{a,b}$,\\
\small $^a$School of Mathematical Sciences\\
\small University of Science and Technology of China, Hefei, Anhui 230026, China.\\
\small $^{b}$ CAS Key Laboratory of Wu Wen-Tsun Mathematics\\
\small University of Science and Technology of China, Hefei, Anhui 230026, China.\\
\small xmhou@ustc.edu.cn
}
\date{}
\maketitle

\begin{abstract}
Given an irreducible fraction $\frac{c}{d} \in [0,1]$, a pair $(\mathcal{A},\mathcal{B})$ is called a $\frac{c}{d}$-cross-intersecting pair of $2^{[n]}$ if $\mathcal{A}, \mathcal{B}$ are two families of subsets of $[n]$ such that for every pair $A \in\mathcal{A}$ and $B\in\mathcal{B}$, $|A \cap B|= \frac{c}{d}|B|$. Mathew, Ray, and Srivastava [{\it\small Fractional cross intersecting families, Graphs and Comb., 2019}] proved that $|\mathcal{A}||\mathcal{B}|\le 2^n$ if $(\mathcal{A}, \mathcal{B})$ is a $\frac{c}{d}$-cross-intersecting pair of $2^{[n]}$ and characterized all the pairs $(\mathcal{A},\mathcal{B})$ with $|\mathcal{A}||\mathcal{B}|=2^n$, such a pair also is called a maximal $\frac cd$-cross-intersecting pair of $2^{[n]}$, when $\frac cd\in\{0,\frac12, 1\}$. In this note, we characterize all the maximal $\frac cd$-cross-intersecting pairs $(\mathcal{A},\mathcal{B})$ when $0<\frac{c}{d}<1$ and $\frac cd\not=\frac 12$, this result answers a question proposed by Mathew, Ray, and Srivastava (2019).

\end{abstract}

\section{Introduction}
We write $[n]$ for $\{1,2,...,n\}$ and $2^{[n]}$ for the power set of $[n]$.  An $\mathcal{F} \subset 2^{[n]}$ is called an {\em intersecting family} if the intersection of every two sets in $\mathcal{F}$ is non-empty. Let $\binom {[n]}{k}$ be the family of all $k$-subsets of $[n]$. An intersecting family $\mathcal{F}$ is called $k$-uniform if $\mathcal {F}\subset \binom {[n]}{k}$.
The famous Erd\H{o}s-Ko-Rado Theorem  states that:
\begin{thm}[Erd\H{o}s-Ko-Rado Theorem, \cite{EKR}]\label{THM: EKR}
$|\mathcal {F} |\leq \binom{n-1}{k-1}$ if $\mathcal{F}$ is a $k$-uniform intersecting family for $n\geq 2k$. Moreover,  the equality holds when $n>2k$ if and only if $\mathcal{F}=\{F\in \binom {[n]}{k} : i\in F\}$ for some $i\in[n]$.
\end{thm}
There are several extensions of Theorem~\ref{THM: EKR} in literatures. Let $L=\{\ell_1,\ell_2,...,\ell_s\}$ be a set of $s$ non-negative integers. A family $\mathcal{F}\subset 2^{[n]}$ is called $L$-intersecting if for every pair of different sets $F_i, F_j\in \mathcal{F}$, $|F_i\cap F_j| \in L$. Let $\mathcal{F}$ be a $t$-uniform $L$-intersecting family. If $L$ is singleton, Bose~\cite{B} showed that $|\mathcal{F}|\leq n$, and in general, Ray-Chaudhuri and Wilson~\cite{RW} showed that $|\mathcal{F}|\leq \binom{n}{s}$, and the upper bound can be achieved by the family $\mathcal{F}=\binom{[n]}{s}$ with $L=\{0,1,...,s-1\}$. Frankl and Wilson~\cite{FW} extended the above result to the non-uniform case by showing that $|\mathcal{F}|\leq \binom{n}{0}+ \binom{n}{1}...+\binom{n}{s}$. More extensions can be found in~\cite{ABS91,BMM, GS02, LY, QR-C00, S94, S03}.
%In %2019, Balachadran, Mathew and Mishra introduced a fractional variant of the classic families in \cite{BMM}.
%See \cite{LY} for a survey on intersecting families.

In this note, we concern another variant of intersecting families: the cross-intersecting families. Two families $\mathcal{A}, \mathcal{B} \subset 2^{[n]}$ is {\em cross-intersecting} if  $A \cap B \neq \emptyset$ for every pair $A\in \mathcal{A}, B\in \mathcal{B}$, $(\mathcal{A},\mathcal{B})$ is also called a {\em cross-intersecting pair}. A cross-intersecting version of the Erd\H{o}s-Ko-Rado Theorem was first given by Pyber~\cite{P}, and   Frankl et al.~\cite{FLST} gave a generalized version by showing that if $\mathcal{A},\mathcal{B} \subset \binom{[n]}{k}$ such that $|A\cap B| \geq t$ for all $A\in \mathcal{A}, B\in \mathcal{B}$, then for all $n\geq (t+1)(k-t+1)$, $|\mathcal{A}||\mathcal{B}|\leq {\binom{n-t}{k-t}}^2$.
A cross-intersecting pair $(\mathcal{A},\mathcal{B})$ with $\mathcal{A},\mathcal{B}\subset 2^{[n]}$ is said to be $\ell$-cross-intersecting if $|A\cap B|=\ell$ for some positive integer $\ell$ and all of  $A\in \mathcal{A}$ and $B\in \mathcal{B}$. The $\ell$-cross-intersecting families were also studied in literatures, for example, in~\cite{ACZ89,AL09}.  Recently, Mathew, Ray, Srivastava~\cite{RRS} introduced a fractional variant of the cross-intersecting family.
Let $\frac{c}{d} \in [0,1]$ be an irreducible fraction. We call $(\mathcal{A},\mathcal{B})$  a $\frac{c}{d}$-cross-intersecting pair of $2^{[n]}$ if $\mathcal{A}, \mathcal{B}\subseteq 2^{[n]}$ and for every pair $A\in \mathcal{A}$ and $B \in \mathcal{B}$, $|A \cap B|= \frac{c}{d}|B|$.
%We call such families pair $({\mathcal{A},\mathcal{B}})$ $\displaystyle\frac{c}{d}$-intersecting families pair.
Let
$${\mathcal{M}_{\frac{c}{d}}}(n)=\max\{|\mathcal{A}||\mathcal{B}| : (\mathcal{A},\mathcal{B}) \mbox{ is a $\frac{c}{d}$-cross-intersecting pair of $2^{[n]}$}\}.$$
Clearly, $\mathcal{M}_{\frac{c}{d}}(n)\ge 2^n$ because $(\mathcal{A},\mathcal{B})$ with $\mathcal{A}=2^{[n]}$ and $\mathcal{B}=\{\emptyset\}$ is a trivial $\frac{c}d$-cross-intersecting pair with $|\mathcal{A}||\mathcal{B}|=2^n$.
We call a $\frac{c}{d}$-cross-intersecting pair $(\mathcal{A},\mathcal{B})$ of $2^{[n]}$ a {\em maximal pair} if $|\mathcal{A}||\mathcal{B}|=\mathcal{M}_{\frac{c}{d}}(n)$.
Mathew, Ray, Srivastava~\cite{RRS} proved that the lower bound is the exact value of $\mathcal{M}_{\frac{c}{d}}(n)$.
\begin{thm}[\cite{RRS}]\label{THM:m_c/d}
For any given irreducible fraction $\frac{c}d\in[0,1]$ and positive integer $n$, $${\mathcal{M}_{\frac{c}{d}}}(n)=2^{\displaystyle n}.$$
\end{thm}
Moreover, they characterized all maximal pairs  when $\frac{c}{d}\in\{0,1,\frac{1}{2}\}$.
%We call $(\mathcal{A},\mathcal{B})$ with $\mathcal{A}=2^{[n]}$ and $\mathcal{B}=\{\emptyset\}$ a trivial maximal pair.
\begin{thm}[\cite{RRS}]\label{THM: c/d=0,1,1/2}
Let $(\mathcal{A},\mathcal{B})$ be a maximal $\frac{c}{d}$-cross-intersecting pair of $2^{[n]}$. Then the following holds.

(1) If $\frac{c}d=0$ then $(\mathcal{A},\mathcal{B})=(2^{[k]}, \mathcal{P}(S))$ for some $0\le k\le n$, where $\mathcal{P}(S)$ is the power set of $\{k+1, \ldots, n\}$.

(2) If $\frac{c}d=1$ then $(\mathcal{A},\mathcal{B})=\left([k]\cup T, 2^{[k]}\right)$ for some $0\le k\le n$,  where $T\in\mathcal{P}(S)$ and $S=\{k+1, \ldots, n\}$.

(3) If $\frac cd=\frac12$, then $(\mathcal{A},\mathcal{B})$ is one of the following $\lfloor\frac{n}{2}\rfloor+1$ pairs of families $({\mathcal{A}}_k,{\mathcal{B}}_k)$, $0\leq k \leq \lfloor\frac{n}{2}\rfloor$,

$\mathcal{A}_0=2^{[n]}$ and $\mathcal{B}_0=\emptyset$,

$\mathcal{A}_k=\{A\in 2^{[n]} : |A\cap\{2i-1,2i\}|=1, \mbox{ for all } i\in [k]\}$,

$\mathcal{B}_k=\{B \in 2^{[n]} : |B\cap \{2i-1,2i\}|\in\{0,2\}, \mbox{ for all } i\in [k]   \mbox{ and for all } j\geq 2k, j\notin B\}$.

The structures of $(\mathcal{A},\mathcal{B})$ are unique, up to isomorphism.
\end{thm}

The authors in~\cite{RRS} also proposed the following natural and interesting problem.
\begin{prob}[Mathew, Ray, Srivastava~\cite{RRS}]\label{PROB: p1}
It would be interesting to show a characterization theorem for any $0<\frac{c}{d}<1$ and $\frac{c}d\not=\frac 12$.
\end{prob}
In this note, we solve the above problem.
%characterize the maximal pairs when $\displaystyle\frac{c}{d}\neq \displaystyle\frac{1}{2}$.
\begin{thm}\label{THM: main}
Suppose $0<\frac{c}{d}<1$ and $\frac{c}d\not=\frac 12$. Then
$(\mathcal{A},\mathcal{B})$ is a maximal $\frac{c}{d}$-cross intersecting pair of $2^{[n]}$ if and only if $(\mathcal{A}, \mathcal{B})=(2^{[n]}, \emptyset)$.
\end{thm}

The rest of the note is arranged as follows. We give some notation and preliminaries in Section 2. The proof of Theorem~\ref{THM: main} will be given in Section 3. We finish this work with some discussion and remarks.

%Together with the theorem,we give a characterization for all maximal pairs.

%\begin{thm}\label{thm4}
%  Let $(\mathcal{A},\mathcal{B})$ be a $\frac{c}{d}$-cross intersecting pair of families of subsets of $[n]$ with $|\mathcal{A}||\mathcal{B}|=2^n$. If $\frac{c}{d} \neq \frac{1}{2}$ then $\mathcal{A}=2^{[n]}$ and $\mathcal{B}=\emptyset$

%If $\frac{c}{d} = \frac{1}{2}$, then $(\mathcal{A},\mathcal{B})$ is one of the following $\lfloor {\frac{n}{2}+1} \rfloor$ pairs of families $({\mathcal{A}}_k,{\mathcal{B}}_k)$,$0\leq k \leq \lfloor {\frac{n}{2}+1} \rfloor$, up to isomorphism,

%$\mathcal{A}_0=2^{[n]}$ and $\mathcal{B}_0=\emptyset$

%$\mathcal{A}_k=\{ \mathcal{A}_k \in 2^{[n]}:|A\cap\{2i-1,2i\}|=1,\forall i, 1\leq i \leq k\}$

%$\mathcal{B}_k=\{ \mathcal{B}_k \in 2^{[n]}: |B\cap \{2i-1,2i\} \in \{0,2\} \forall i, 1\leq i\leq k \quad and \quad \forall j\geq 2k, j\notin B$
%\end{thm}

\section{Preliminaries}
For any $S \subseteq [n]$, let $X_S\in\mathbb{R}^n$ denote the characteristic vector of $S$ and $X_S(i)$ denote its $i$-th entry, i.e.
$$ X_S(i)=\left\{
\begin{array}{rcl}
0      &     & {i\notin S}\\
1    &     & {i\in    S}\\
\end{array} \right.. $$
The weight of a characteristic vector is the number of its non-zero entries. So the weight of $X_S$ is equal to $|S|$.
For any family $\mathcal{A}\subseteq 2^{[n]}$, we shall not distinguish the family and the collection of their corresponding characteristic vectors. %This meaning will be clearly stated from the context.

Given $V\subseteq \mathbb{F}_2^n$, let $\span(V)$ be the vector subspace spanned by $V$.
%the all linear combinations of the vectors of $V$ in $\mathbb F^{2}_{n}$.
We will use $\basis(V)$ denote a basis of $\span(V)$ and $\dim(V)$ the dimension of $\span(V)$.
A collection of vectors $V\subseteq\mathbb{F}_2^n$ is called a {\it linear code} if $V=\span(V)$.
Write `$\langle\  \rangle_{1}$' and `$\langle\  \rangle_2$' for the inner products in $\mathbb{R}^n$ and in  $\mathbb{F}_2^n$, respectively, i.e.
$$\langle x,y\rangle_{1}=x_1y_1+x_2y_2+...+x_ny_n \mbox{ for } x,y\in\mathbb{R}^n,$$
and
$$\langle x,y\rangle_{2}=x_1y_1+x_2y_2+...+x_ny_n\pmod2 \mbox{ for } x,y\in\mathbb{F}_2^n.$$
Given a linear code $C\subseteq \mathbb F^{2}_{n}$, the {\it dual code} $C^{\perp}$ is defined as
$$C^{\perp}=\{x\in \mathbb F^{2}_{n} \, :\, \langle x,c\rangle_2=0 \text{ for all } c\in C\}.$$
%where $<x,y>_2$ is the standard inner product in $\mathbb{F}_2^n$, i.e. for $x=(x_1,x_2,...,x_n),y=(y_1,y_2,...,y_n)$,
%$$<x,y>_2=x_1y_1+x_2y_2+...+x_ny_n      (mod \quad 2)$$.
It is easy to verify the following property.
\begin{prop}[\cite{RRS}]\label{PROP: dualC}
 If $C$ is a linear code, then $C^{\perp}$ is a linear code too.
\end{prop}

We need a characterization when the binomial coefficient is a power of 2.
\begin{prop}\label{Prop: nchoosek}
$\binom{n}{k}$ is a power of $2$ if and only if $k=0$ or $(k,n)\in\{(1,2^m), (2^m-1, 2^m)\}$.
\end{prop}
\begin{proof}
%The proof is straightforward.We need only to use the prime decompositon to check the divisor $2$ in the combination number.
Suppose $\binom{n}{k}=\frac{n!}{k!(n-k)!}=2^m$ for some integer $m\ge 0$. So
\begin{eqnarray*}
m&=&\sum_{2^i\leq n}\left\lfloor\frac{n}{2^i}\right\rfloor-\sum_{2^i\leq k}\left\lfloor\frac{k}{2^i}\right\rfloor-\sum_{2^i\leq n-k} \left\lfloor\frac{n-k}{2^i}\right\rfloor\\
&=&\sum_{2^i\leq n}\left(\left\lfloor\frac{n}{2^i}\right\rfloor-\left\lfloor\frac{k}{2^i}\right\rfloor-\left\lfloor\frac{n-k}{2^i}\right\rfloor\right)\\
&\leq& \lfloor\log_2 n\rfloor,
\end{eqnarray*}
the last inequality holds because $\left\lfloor\frac{n}{2^i}\right\rfloor-\left\lfloor\frac{k}{2^i}\right\rfloor-\left\lfloor\frac{n-k}{2^i}\right\rfloor\le 1$.
Therefore, $\binom{n}{k}=2^m$ if and only if $2^m=\binom{n}{k}\leq 2^{\log_2 n}=n$ if and only if $m=k=0$, or $(k,n)=(1,2^m)$, or $(k,n)=(2^m-1, 2^m)$.
\end{proof}

The following lemma has been proved in~\cite{RRS}.
\begin{lem}[\cite{RRS}]\label{LEM: matrixM}
If the elments of a linear code $C\subset F_2^n$ are arranged as rows of a matrix $M_C$ with $n$ columns, then for each column, one of the following holds.
%\begin{enumerate}

\noindent(1) All the entries in that column are 0.

\noindent(2) Exactly half the entries in that column are 0, and the rest are 1.
%\end{enumerate}
\end{lem}

\section{Proof of Theorem~\ref{THM: main}}

Suppose $(\mathcal{A},\mathcal{B})$ is a maximal $\frac{c}{d}$-cross-intersecting pair of $2^{[n]}$ and $(\mathcal{A},\mathcal{B})\not=(2^{[n]}, \emptyset)$, where $\frac cd$ is an irreducible fraction with $0<\frac{c}{d}<1$ and $\frac cd\not=\frac 12$. Then, by Theorem~\ref{THM:m_c/d}, we have $|\mathcal{A}||\mathcal{B}|=2^n$.
%is an irreducible fraction.a maximal pair.
As shown in the proof of Theorem~\ref{THM:m_c/d} in~\cite{RRS}, we partition $\mathcal{B}$ into two parts
$$\mathcal{B}_1 = \{B\in \mathcal{B} : |B|\equiv 0\pmod {2d}\},$$
$$\mathcal{B}_2 = \{B\in \mathcal{B} : |B|\equiv d\pmod {2d}\}.$$
Recall that  $\langle X_A, X_B\rangle_{1}=|A\cap B|=\frac cd|B|$ is an integer for any $A\in \mathcal{A}, B\in \mathcal {B}$, we get
%$<X_A, X_B>_{2}=1$ if $c$ is odd and $B\in \mathcal{B}_2$; $<X_A, X_B>_{2}=0$,otherwise
\begin{equation}\label{EQ: e1}
\langle X_A, X_B\rangle_{2}=\left\{
\begin{array}{rcl}
1      &     & {B\in\mathcal{B}_2 \text{ and $c$ is odd}}\\
0    &     & \text{otherwise}\\
\end{array} \right..
\end{equation}

If $c$ is even, then $\langle X_A, X_B\rangle_{2}=0$ for all $A\in \mathcal{A}, B\in \mathcal{B}$. Thus $\span(\mathcal{A}) \bot \span(\mathcal{B})$ in $\mathbb{F}_2^n$ and so $\dim(\span(\mathcal{A}))+\dim(\span(\mathcal{B})) \leq n$. Therefore,
$$2^n=|\mathcal{A}| |\mathcal{B}| \leq |\span(\mathcal{A})||\span(\mathcal{B})|=2^{\dim(\span(\mathcal{A}))}2^{\dim(span(\mathcal{B}))}=2^{\dim(\span(\mathcal{A}))+\dim(\span(\mathcal{B}))}\leq 2^{n}.$$
So we have
\begin{equation}\label{EQU: mathcal{A}}
\span{\mathcal{(A)}}=\mathcal{A},\,  \span{\mathcal{(B)}}=\mathcal{B}, \text{ and }\dim(\span(\mathcal{A}))+\dim(\span(\mathcal{B}))=n.
\end{equation}

If $c$ is odd, we construct $\mathcal{B}_1'$ by appending a 0 to the left of every vector in $\mathcal{B}_1$, and $\mathcal{B}_2'$ by appending a 1 to the left of every vector in $\mathcal{B}_2$. Let $\mathcal{B}'=\mathcal{B}_1' \cup \mathcal{B}_2'$. Construct $\mathcal{A}'$ by appending 1 to the left of every vector in $\mathcal{A}$.
Then $\langle X_A, X_B\rangle_2=0$ for all $A\in \mathcal{A}', B\in \mathcal{B}'$.
So $\span(\mathcal{A}')\bot\span(\mathcal{B}')$ in $\mathbb {F}_2^{n+1}$.
Therefore, we have $$\dim(\span{\mathcal{(A')}})+\dim(\span{\mathcal{(B')}})\leq n+1$$
and
$$|\span(\mathcal{A}')||\span(\mathcal{B}')|=2^{\dim(\span(\mathcal{A}'))}2^{\dim(span(\mathcal{B}'))}
=2^{\dim(\span({\mathcal{A}'}))+\dim(\span({\mathcal{B}'}))}\leq 2^{n+1}.$$

\begin{claim}\label{CLM: claim1}
$|\span({\mathcal{A}'})|=2|\mathcal{A}'|$, $\span({\mathcal{B}'})=\mathcal{B}'$ and $\dim(\span({\mathcal{A}'}))+\dim(\span({\mathcal{B}'}))=n+1$.
\end{claim}
\begin{proof}
We first claim that $|\span({\mathcal{A}'})|\geq 2|{\mathcal{A}'}|$. In fact, as $\span({\mathcal{A}'})$ is a linear code, by Lemma~\ref{LEM: matrixM}, the leftmost column of the matrix $M_{\span(\mathcal{A}')}$ does contain $0$. But the leftmost entry of all the vectors in $\mathcal{A}'\subseteq\span(\mathcal{A}')$ are  $1$, so there are at least $|\mathcal{A}'|$ vectors in $\span(\mathcal{A}')$ having their leftmost entry as $0$, i.e. $|\span({\mathcal{A}'})|\geq 2|{\mathcal{A}'}|$.
Therefore,
$$2^{n+1}=2|\mathcal{A}||\mathcal{B}|=2|\mathcal{A}'||\mathcal{B}'|\le |\span({\mathcal{A}'})||\span({\mathcal{B}'})|\leq 2^{n+1}.$$
This implies the claim.
\end{proof}

Let $f_{\mathcal{A}}: \mathcal{A}\mapsto \mathcal{A}'$ be the bijection that maps every vector in $\mathcal{A}$ to the corresponding vector in $\mathcal{A}'$, and let $g_{\mathcal{A}}$ be its inverse. Similarly, we can define $f_{\mathcal{B}}$ and its inverse $g_{\mathcal{B}}$, i.e.,
%For any subset $V\subset \mathcal{A}$, we use $f_\mathcal{A}(V)$ to denote $\{f_\mathcal{A}(A)| A\in V\}$. For any subset $V'\subset \mathcal{A}'$,we use $g_\mathcal{A}$ to denote set .
%We can define $f_\mathcal{B}(V)$ and $g_\mathcal{B}(V')$ similarly.
we have $f_\mathcal{B}(\mathcal{B}_1)=\mathcal{B}_1'$ and $f_\mathcal{B}(\mathcal{B}_2)=\mathcal{B}_2'$; $g_\mathcal{B}(\mathcal{B}_1')=\mathcal{B}_1$ and $g_\mathcal{B}(\mathcal{B}_2')=\mathcal{B}_2$.

\begin{claim}\label{CLM:claim2}
$\mathcal{B}$ is a linear subspace of $\mathbb{F}_2^n$.
\end{claim}
\begin{proof}
If $c$ is even, the result follows directly from (\ref{EQU: mathcal{A}}).

Now suppose $c$ is odd. Then, by Claim~\ref{CLM: claim1}, $\mathcal{B}'=\span(\mathcal{B}')$ is a linear space. Let $B_1, B_2\in \mathcal{B}$. Note that $X_{B_1}+X_{B_2}=X_{B_1\Delta B_2}$ in $\mathbb{F}_2^n$. Denote $B_3=B_1\Delta B_2$.  It is sufficient to show that $B_3 \in \mathcal{B}$.
%This equality shows the $\mathcal{B}$ is closed under addition in $F2n$ over $F2$, hence $|span(\mathcal{B})|=|\mathcal{B}|$.
Let $B_1'=f_\mathcal{B}(B_1), B_2'=f_\mathcal{B}(B_2)$. Then $B_3'=B_1'\Delta B_2'\in \mathcal{B}'$. Obviously, $X_{B_3'}$ is obtained by appending a $0$ (if $B_1, B_2$ are in a same part $\mathcal{B}_i$ for some $i=1,2$) or $1$ (otherwise) to the leftmost entry of $X_{B_3}$, i.e. $B_3=g_\mathcal{B}(B_3')\in\mathcal{B}$.
%Since $x_3\in \mathcal{B}'$, we have $x_3=g_\mathcal{B}(x_3')\in \mathcal{B}$.So $\mathcal{B}$ is a linear space and $|span(\mathcal{B})|=|\mathcal{B}|$.
%For $B_1,B_2$ symetric deffrence $B_3=B_1\Delta B_2$, $b_3=b_1+b_2 \in \mathcal{B}$, so $B_3 \in \mathcal{B}$.
\end{proof}

\begin{claim}\label{CLM: $b_1+b_2$}
For any  $B_1, B_2\in \mathcal{B}$,

(1) $
B_1\Delta B_2\in\left\{\begin{array}{ll}
\mathcal{B}_1 &\text{ either $B_1, B_2\in \mathcal{B}_1$, or $B_1, B_2 \in \mathcal{B}_2$}\\
\mathcal{B}_2 & \text{ otherwise}
\end{array}\right.
$;

(2) $|B_1\cap B_2|\equiv 0 \pmod d$.
\end{claim}
%\proof
%The forward direction is straightforword by the construction. For the opposite direction,
%$b_1+b_2\in B_1$:  $f_\mathcal{B}(b_1+b_2)=f_\mathcal{B}(b_1)+f_\mathcal{B}(b_2)=b_1'+b_2'\in B_1'$. By the construction, the lefemost bits of $b_1'$ and $b_2'$ are both $0$ or $1$, this means that $b_1$ and $b_2$ are both from $\mathcal{B}_1$ or $\mathcal{B}_2$.
%\begin{thm}
%$|B_1\cap B_2|=0 (mod \quad d)$
%\end{thm}
\begin{proof}
 (1) follows directly from the proof of Claim~\ref{CLM:claim2}.

 (2) It follows from (1) and the fact that $|B_1\Delta B_2|=|B_1|+|B_2|-2|B_1\cap B_2|$.
 %From the proposition, we have $B_1\Delta B_2 \in \mathcal{B}$. When $B_3 \in \mathcal{B}_1$,then both $B_1$ and $B_2$ either in $\mathcal{B}_1$ or in $\mathcal{B}_2$.Therefore, $|B_1\cap B_2|= 0 (mod\quad d)$.
%When $B_3 \in \mathcal{B}_2$, then either $B_1\in \mathcal{B}_1, B_2\in \mathcal{B}_2$ or $B_1\in \mathcal{B}_2,B_2\in \mathcal{B}_1$.Therefore, $|B_1\cap B_2|=0 (mod\quad d)$.
\end{proof}

\begin{claim}\label{CLM：B1capB2}
%(1) If $A\in\mathcal{A}$ and $B_1, B_2\mathcal{B}$$|A\cap B_1|=\frac{c}{d}|B_1|$, $|A\cap B_2|=\frac{c}{d}|B_2|$ and $|A\cap (B_1\Delta B_2)|=\frac{c}{d}|B_1\Delta B_2|$, then $|A\cap (B_1\cap B_2)|=\frac{c}{d}|B_1\cap B_2|$.
$\mathcal{B}$ is closed under intersection.
\end{claim}
%\proof
%$|A\cap B_1\Delta B_2|=|A\cap((B_1 b_2)\cup (B_2 B_1))|=\frac{c}{d}|B_1\Delta B_2|$
%$|A\cap (B_1 B_2)|+|A\cap (B_2 B_1)|=\frac{c}{d}(|B_1|+|B_2|-2|B_1\cap B_2|)$
%$|A\cap B_1|-|A\cap (B_1\cap B_2)|+|A\cap B_2|-|A\cap (B_1\cap B_2)|=\frac{c}{d}(|B_1|+|B_2|-2|B_1\cap B_2|$
%$\frac{c}{d}|B_1|+\frac{c}{d}|B_2|-2|A\cap (B_1\cap B_2)|=\frac{c}{d}(|B_1|+|B_2|-2|B_1\cap B_2|$
%So we have $|A\cap |B_1\cap B_2||=\frac{c}{d}|B_1\cap B_2|$.
%\begin{prop}
%$\mathcal{B}$ is closed under intersection.
%\end{prop}
\begin{proof}
It is sufficient to show that $B_1\cap B_2 \in \mathcal{B}$ for any $B_1, B_2\in \mathcal{B}$. By Claim~\ref{CLM: b_1+b_2}, $B_1\Delta B_2 \in \mathcal{B}$. Hence, for any $A\in\mathcal{A}$, we have $|A\cap (B_1\Delta B_2)|=\frac cd|B_1\Delta B_2|$. So
$$2|A\cap (B_1\cap B_2)|=|A\cap B_1|+|A\cap B_2|-|A\cap (B_1\Delta B_2)|=\frac{c}{d}(|B_1|+| B_2|-|B_1\Delta B_2|)=2\cdot\frac cd|B_1\cap B_2|.$$ Since $(\mathcal{A},\mathcal{B})$ is a maximal pair, $B_1\cap B_2 \in \mathcal{B}$.
\end{proof}

A set $B\in\mathcal{B}$ is called {\it primitive} if for any set $B'\in\mathcal{B}$, $B'\cap B=\emptyset$ or $B$. Obviously, two different primitive sets are disjoint. Let  $\mathcal{S}=\{B_1,B_2,...,B_k\}$ be the set of all nonempty primitive sets in $\mathcal{B}$.
%With these propsition, it's enough to characterize $\mathcal{B}$
\begin{claim}\label{CLM: basis}
$\mathcal{S}=\{B_1,...,B_k\}$ is a basis of $\mathcal{B}$.
%There exist a basis $\mathcal{S}=\{B_1,...,B_k\}\subset \mathcal{B}$ such that (i) the set is $\mathcal{S}$ are pairwise disjoint, and (ii) for every set $B\in \mathcal{B}$, there exist some sets in $\mathcal{S}$ whose union yields $B$.
\end{claim}
\begin{proof}
%From the prop, we see that $B$ is closed under intersection and symmetric difference. Define a set $B$ is primitive if for any set in $\mathcal{B}$ the intersection with $B$ is either empty or $B$. Obviously two different primitive sets are disjoint.Let  $\mathcal{S}=\{B_1,B_2,...,B_k\}$ be the set of all primitive sets.
We first claim that every element in any set $B\in\mathcal{B}$ must be present in exactly one set in $\mathcal{S}$.
Obviously, an element in $B$ cannot be present in more than one set in $\mathcal{S}$ for two primitive sets are disjoint. Now we suppose that there exist some elements in $B$ that are not present in any set in $\mathcal{S}$.  Since $\mathcal{B}$ is closed under intersection and $B_1\cup \ldots \cup B_k=B_1\Delta \ldots\Delta B_k\in\mathcal{B}$, there must exist nonemptysets (e.g. $B\cap (B\Delta (B_1\cup\ldots\cup B_k))$) with no element present in any set in $\mathcal{S}$. Choose a smallest one, say $B_{\min}$, of  such sets in  $\mathcal{B}$. By the definition of $\mathcal{S}$, $B_{\min}$ can not be primitive.
%Then $B_{min}\cap B_i=\emptyset$ for every $B_i\in \mathcal{B}$. If not, $B_{min}\Delta B_{i}=B_{min}\backslash B_i$  is a set smaller than $B_{min}$, a contradiction. Thus, $B_{min}$ is disjoint from all primitive sets $B_i$ in $\mathcal{B}$. But the set $B_{min}$ is not primitive,
So there must exist a set $B'\in \mathcal{B}$ such that $B_{\min}\cap B'$ is neither $\emptyset$ nor $B_{\min}$. Clearly, $B_{\min}\cap B'$ has no element present in any set in $\mathcal{S}$, but $|B_{\min}\cap B'|<|B_{\min}|$, which is a contradiction.
This claim also implies that $B\subseteq B_1\cup B_2\cup\ldots\cup B_k$ for any set $B\in\mathcal{B}$. To show $\mathcal{S}$ is a basis of $\mathcal{B}$, it is sufficient to prove that $B=B_{i_1}\Delta\ldots\Delta B_{i_\ell}$ for some $i_1,\ldots,i_\ell\in\{1,2,\ldots,k\}$.
In fact, let $B_{i_1},\ldots, B_{i_\ell}$ be all sets with $B_{i_j}\cap B\not=\emptyset$. Then $B_{i_j}\subseteq B$ for any $1\le j\le \ell$ as, otherwise, $B\cap B_{i_j}\neq B_{i_j}$ and $\emptyset$, which is a contradict to the primitivity of $B_{i_j}$.
%  Up to isomorphism, let there be $n_0$ elements, say $\{n-n_0+1,...,n\}$ are not present in any set in $\mathcal{B}$. The remaining $n-n_0$ elements ($[n-n_0]$) appear in exactly one set in $\mathcal{S}$, so we can regard $\mathcal{S}$ the partition of $[n-n_0]$. For any set $B$ in $\mathcal{B}$, $B$ is a subet of $[n-n_0]$ and it is able to be written as a series of primitive sets. If not,
\end{proof}

By Claim~\ref{CLM: basis}, we have $|\mathcal{B}|=2^k$. Up to isomorphism, we may assume $B_1\cup B_2\cup\ldots\cup B_k=[n-n_0]$ for some integer $n_0$.
%let there be $n_0$ elements, say $\{n-n_0+1,...,n\}$ are not present in any set in $\mathcal{B}$. The remaining $n-n_0$ elements ($[n-n_0]$) appear in exactly one set in $\mathcal{S}$,
By the definition of $\mathcal{S}$, $B_1, B_2, \ldots, B_k$ is a partition of $[n-n_0]$.
%For any set $B$ in $\mathcal{B}$, $B$ is a subet of $[n-n_0]$ and it is able to be written as a series of primitive sets. If not,With the theorem, we have a special basis $\mathcal{S}=\{B_1,...,B_k\}$ which the elements are disjoint of $\mathcal{B}$. It's a partition of $[n-n_0]$.
Since $(\mathcal{A},\mathcal{B})$ is a $\frac{c}{d}$-cross intersecting pair, the size of every set in $\mathcal{B}$ is divided by $d$. Suppose $|B_1|=d\ell_1$, $|B_2|=d\ell_2$, \ldots, $|B_k|=d\ell_k$. Then $n-n_0=d(\ell_1+...+\ell_k)$. As for every $A\in\mathcal{A}$, we have $|A\cap B_i|=c\ell_i$ for some $0\le c\le \ell_i$ and $A\cap\{n-n_0+1, n-n_0+2,\ldots, n\}$ can be chosen arbitrarily from $\{n-n_0+1, n-n_0+2,\ldots, n\}$.
So $|\mathcal{A}|=2^{n_0}\cdot\prod_{i=1}^k\binom{d\ell_i}{c\ell_i}$. Therefore, we have
%\begin{claim}\label{CLM: AcapB_i}
%$|A\cap B|=\frac{c}{d}|B|$ if and only if $|A\cap B_i|=\frac{c}{d}|B_i|$.
%\end{claim}
%\proof
%The proof of the necessary condition is straightforward for $\mathcal{S} \subset \mathcal{B}$. For the sufficient condition, let $B=B_1\cup B_2\cup ...\cup B_l$.
%We have $|A\cap B|=|A\cap(\cup B_j)|=\Sigma |A\cap B_j|=\Sigma \frac{c}{d}|B_j|=\frac{c}{d}|\cup B_j|=\frac{c}{d}B$.
%Now we can give a characterization. Since each set in $\mathcal{A}$ are disjoint, and $|B_1|=dl_1$,$|B_2|=dl_2$,...,$|B_k|=dl_k$, so there are $2^{n_0}$ sets in $\mathcal{A}$. By the claim 2, we can calculate
\begin{equation}\label{eq1}
2^n=|\mathcal{A}||\mathcal{B}|=2^{n_0}\cdot\prod_{i=1}^k\binom{d\ell_i}{c\ell_i}\cdot 2^k.
\end{equation}
%From theorem, we know the bound is $2^n$. $2^n$ is divided by 2(prime), so
The Equality (\ref{eq1}) implies that every combinatorial number $\binom{d\ell_i}{c\ell_i}$ in the right hand must be a power of $2$.
Note that $\frac cd$ is an irreducible fraction with $0<\frac{c}{d}<1$ and $\frac cd\not=\frac 12$.
If there is some $i\in[k]$ with $\binom{d\ell_i}{c\ell_i}=2^0=1$, then $d\ell_i=c\ell_i$. So $\ell_i=0$ as $1\le c<d$. This is a contradiction to $B_i\not=\emptyset$. So for all $i\in[k]$, we have $\binom{d\ell_i}{c\ell_i}=2^{m_i}$ for some integer $m_i\ge 1$.
\begin{claim}\label{CLM: d=2^m}
For all $i\in[k]$, we have $\ell_i=1$ and $d=2^m$ for some integer $m\ge 2$.
\end{claim}
\begin{proof}
By Proposition~\ref{Prop: nchoosek}, $(c\ell_i, d\ell_i)=(1, 2^{m_i})$ or $(2^{m_i}-1, 2^{m_i})$. For the former case, we have $\ell_i=c=1$ and so $d=2^{m_i}$. For the latter, we have $c\ell_i+1=d\ell_i=2^{m_i}$ and so $\ell_i=1$ and $c+1=d=2^{m_i}$ as $c<d$.
Since $d$ is a constant, we have $d=2^m$ for some positive integer $m$. $m\ge 2$ because  $\frac cd\not=\frac 12$.
\end{proof}
By Claim~\ref{CLM: d=2^m}, for all $i\in[k]$,  we have $\binom{d\ell_i}{c\ell_i}=2^{m}$ for some integer $m\ge 2$. By Equality (\ref{eq1}),
$$2^n=|\mathcal{A}||\mathcal{B}|=2^{n_0}\cdot\prod_{i=1}^k\binom{d\ell_i}{c\ell_i}\cdot 2^k=2^{n_0+(m+1)k}<2^{n_0+k2^m}=2^n,$$ for all $m\ge 2$, a contradiction.
This completes the proof of Theorem~\ref{THM: main}.

\section{Remarks and Discussions}	
In this note, we characterize the maximal $\frac{c}{d}$-cross-intersecting pairs $(\mathcal{A},\mathcal{B})$ of $2^{[n]}$ when $0<\frac{c}{d}<1$ and $\frac{c}d\not=\frac 12$ (Theorem~\ref{THM: main}), this result answers a question proposed by Mathew, Ray, Srivastava~\cite{RRS}. Combining with the result given by Mathew, Ray, Srivastava~\cite{RRS} (Theorem~\ref{THM: c/d=0,1,1/2}), the problem of characterizing the maximal $\frac{c}{d}$-cross intersecting pairs of $2^{[n]}$ for $\frac cd\in[0,1]$ has been solved completely.

%As we have mentioned in the introduction,
For the further study, we can extend the fractional cross-intersecting families to  a symmetric $\frac ab$-cross-intersecting family as follows.
Let $\mathcal{A}$ and $\mathcal{B}$ be two families of subsets of $[n]$. Given an irreducible fraction $\frac{a}{b} \in [0,1]$, we call $(\mathcal{A}, \mathcal{B})$ a {\it symmetric $\frac ab$-cross-intersecting pair} if for every pair $A\in\mathcal{A}$ and $B\in\mathcal{B}$, $|A \cap B|\in \left\{\frac{a}{b}|A|,\frac{a}{b}|B|\right\}$.
It will be very interesting to determine the maximum of $|\mathcal{A}||\mathcal{B}|$ and characterize the maximal pairs $(\mathcal{A}, \mathcal{B})$ of the symmetric $\frac ab$-cross-intersecting families of $2^{[n]}$. Seemingly this problem is more difficult than the nonsymmetric case studied in this paper.

%The fraction intersecting problem was introduced in \cite{BMM},the intersecting condition is more involved: $\forall F_i,F_j\in\mathcal{F}$,$|F_i\cap F_j|\in \{\frac{c}{d}|F_i|,\frac{c}{d}|F_j|\}$. In this cross intersection problem, the author use the algebra method to build the upper bound. This method heavily depend on the "structure" of the family,i.e.$|A_i \bigcap B_j|= \displaystyle\frac{a}{b}|B_j|$.With the similar condition, we have a such problem:
%\begin{prob}
%Let $\mathcal{A}=\{A_1,A_2,...,A_p\}$ and $\mathcal{B}=\{B_1,B_2,...,B_q\}$ be two families of subsets of $[n]$. For every $i \in [p]$ and $j \in [q]$,$|A_i \bigcap B_j|\in \{\displaystyle\frac{a}{b}|A_i|,\displaystyle\frac{a}{b}|B_j|\}$,where $\displaystyle\frac{a}{b} \in [0,1]$ is an irreducible fraction.What's the upper bound of $|\mathcal{A}||{B}|$?
%\end{prob}
%If we use the same method, it's not able to get a upper bound.It's more involved.

\end{document}